\documentclass[12pt]{article}
\usepackage{amsmath, amsthm, amsfonts}
\usepackage{amssymb}
\usepackage{amscd}
\usepackage{verbatim}

\begin{document}
\newcommand{\dx}{\,\mathrm{d}x}
\newcommand{\core}{C_0^{\infty}(\Omega)}
\newcommand{\sob}{W^{1,p}(\Omega)}
\newcommand{\sobloc}{W^{1,p}_{\mathrm{loc}}(\Omega)}
\newcommand{\merhav}{{\mathcal D}^{1,p}}
\newcommand{\be}{\begin{equation}}
\newcommand{\ee}{\end{equation}}
\newcommand{\mysection}[1]{\section{#1}\setcounter{equation}{0}}
\newcommand{\bea}{\begin{eqnarray}}
\newcommand{\eea}{\end{eqnarray}}
\newcommand{\bean}{\begin{eqnarray*}}
\newcommand{\eean}{\end{eqnarray*}}
\newcommand{\thkl}{\rule[-.5mm]{.3mm}{3mm}}
\newcommand{\cw}{\stackrel{\rightharpoonup}{\rightharpoonup}}
\newcommand{\id}{\operatorname{id}}
\newcommand{\supp}{\operatorname{supp}}
\newcommand{\wlim}{\mbox{ w-lim }}
\newcommand{\mymu}{{x_N^{-p_*}}}
\newcommand{\R}{{\mathbb R}}
\newcommand{\N}{{\mathbb N}}
\newcommand{\Z}{{\mathbb Z}}
\newcommand{\Q}{{\mathbb Q}}
\newcommand{\abs}[1]{\lvert#1\rvert}
\newtheorem{theorem}{Theorem}[section]
\newtheorem{corollary}[theorem]{Corollary}
\newtheorem{lemma}[theorem]{Lemma}
\newtheorem{definition}[theorem]{Definition}
\newtheorem{remark}[theorem]{Remark}
\newtheorem{proposition}[theorem]{Proposition}
\newtheorem{problem}[theorem]{Problem}
\newtheorem{conjecture}[theorem]{Conjecture}
\newtheorem{question}[theorem]{Question}
\newtheorem{example}[theorem]{Example}
\newtheorem{Thm}[theorem]{Theorem}
\newtheorem{Lem}[theorem]{Lemma}
\newtheorem{Pro}[theorem]{Proposition}
\newtheorem{Def}[theorem]{Definition}
\newtheorem{Exa}[theorem]{Example}
\newtheorem{Exs}[theorem]{Examples}
\newtheorem{Rems}[theorem]{Remarks}
\newtheorem{Rem}[theorem]{Remark}
\newtheorem{Cor}[theorem]{Corollary}
\newtheorem{Conj}[theorem]{Conjecture}
\newtheorem{Prob}[theorem]{Problem}
\newtheorem{Ques}[theorem]{Question}
\newcommand{\pf}{\noindent \mbox{{\bf Proof}: }}


\renewcommand{\theequation}{\thesection.\arabic{equation}}
\catcode`@=11
\@addtoreset{equation}{section}
\catcode`@=12


\title{Ground state alternative for $p$-Laplacian with potential term}
\author{Yehuda Pinchover\\
 {\small Department of Mathematics}\\ {\small  Technion - Israel Institute of Technology}\\
 {\small Haifa 32000, Israel}\\
{\small pincho@techunix.technion.ac.il}\\\and Kyril Tintarev
\\{\small Department of Mathematics}\\{\small Uppsala University}\\
{\small SE-751 06 Uppsala, Sweden}\\{\small
kyril.tintarev@math.uu.se}}
\maketitle
\newcommand{\dnorm}[1]{\thkl #1 \thkl\,}

\begin{abstract}
Let $\Omega$ be a domain in $\mathbb{R}^d$, $d\geq 2$, and
$1<p<\infty$. Fix $V\in L_{\mathrm{loc}}^\infty(\Omega)$. Consider
the functional $Q$ and its G\^{a}teaux  derivative $Q^\prime$
given by
$$Q(u):=\!\int_\Omega\!\!\! (|\nabla u|^p+V|u|^p)\!\dx,\;\; \frac{1}{p}Q^\prime
(u):=-\nabla\cdot(|\nabla u|^{p-2}\nabla u)+V|u|^{p-2}u.$$ If
$Q\ge 0$ on $C_0^{\infty}(\Omega)$, then either there is a
positive continuous function $W$ such that $\int
W|u|^p\,\mathrm{d}x\leq Q(u)$ for all $u\in C_0^{\infty}(\Omega)$,
or there is a sequence $u_k\in C_0^{\infty}(\Omega)$ and a
function $v>0$ satisfying $Q^\prime (v)=0$, such that $Q(u_k)\to
0$, and $u_k\to v$ in $L^p_\mathrm{loc}(\Omega$). In the latter
case, $v$  is (up to a multiplicative constant) the unique
positive supersolution of the equation $Q^\prime (u)=0$ in
$\Omega$, and one has for $Q$ an inequality of Poincar\'e type:
there exists a positive continuous function $W$ such that for
every $\psi\in C_0^\infty(\Omega)$ satisfying $\int \psi
v\,\mathrm{d}x \neq 0$ there exists a constant $C>0$ such that
$C^{-1}\int W|u|^p\,\mathrm{d}x\le Q(u)+C\left|\int u
\psi\,\mathrm{d}x\right|^p$ for all $u\in C_0^\infty(\Omega)$. As
a consequence, we prove positivity properties for the quasilinear
operator $Q^\prime$ that are known to hold for general subcritical
resp. critical second-order linear elliptic operators.
\\[2mm]
\end{abstract}

\mysection{Introduction}

Positivity properties of quasilinear elliptic equations, in
particular those with the $p$-Laplacian term in the principal
part, have been extensively studied over the recent decades (see
for example
\cite{AH1,AH2,BFT,BL,BM,DGTU,DKN,FMST,FGT,HKM,MarS,PS,V} and the
references therein). The objective of the present paper is to
study general positivity properties of such equations defined on
general domains in $\R^d$. We generalize some results, obtained
for $p=2$ in \cite{PT2}, to the case of $p\in(1,\infty)$. In
particular, we extend to the case of the $p$-Laplacian the
dichotomy of \cite{PT2} obtained for nonnegative Schr\"odinger
operators which states that either the associated quadratic form
has a weighted spectral gap or the operator admits a unique ground
state.

Some of the proofs in this paper (in particular, the uniqueness of
a global positive supersolution for a critical operator) seem to
be new even for the previously studied case $p=2$.

Fix $p\in(1,\infty)$, a domain $\Omega\subseteq\R^d$, and a
potential $V\in L_\mathrm{loc}^\infty(\Omega)$. We denote by
$\Delta_p(u):=\nabla\cdot(|\nabla u|^{p-2}\nabla u)$ the
$p$-Laplacian.  Throughout this paper we assume that \be \label{Q}
Q(u):=\int_\Omega \left(|\nabla u|^p+V|u|^p\right)\dx\ge 0\ee for
all $u\in \core$.
\begin{definition}{\em We say that a function
$v\in W^{1,p}_{\mathrm{loc}}(\Omega)$ is a {\em (weak) solution}
of the equation
 \be \label{groundstate}
 \frac{1}{p}Q^\prime
(v):=-\Delta_p(v)+V|v|^{p-2}v=0\quad \mbox{in }  \Omega,\ee if
for every $\varphi\in\core$
 \be \label{solution} \int_\Omega (|\nabla v|^{p-2}\nabla
v\cdot\nabla\varphi+V|v|^{p-2}v\varphi)\dx=0. \ee

We say that a positive function $v\in C^1_{\mathrm{loc}}(\Omega)$
is a {\em positive supersolution} of the equation
(\ref{groundstate}) if for every nonnegative $\varphi\in\core$
 \be\label{supersolution}
\int_\Omega (|\nabla v|^{p-2}\nabla
v\cdot\nabla\varphi+V|v|^{p-2}v\varphi)\dx\geq 0. \ee
 }\end{definition}
\begin{remark}{\em It is well-known that any weak solution of (\ref{groundstate})
admits H\"older continuous first derivatives, and that any
nonnegative solution of (\ref{groundstate}) satisfies the Harnack
inequality \cite{DiB,S,T}.
 }\end{remark}

\begin{definition}{\em We say that the functional $Q$ has a
{\em weighted spectral gap} in $\Omega$ (or $Q$ is {\em strictly
positive} in $\Omega$) if there is a positive continuous function
$W$ in $\Omega$ such that  \be\label{wsg} Q(u)\ge \int_\Omega
W|u|^p\dx \qquad \forall u\in\core.\ee }\end{definition}
\begin{definition}{\em We say that a sequence
$\{u_k\}\subset\core$  is a {\em null sequence}, if $u_k\ge 0$ for
all $k\in\N$, and there exists an open set $B\Subset\Omega$ (i.e.,
$\overline{B}$ is compact in $\Omega$) such that
$\int_B|u_k|^p\dx=1$, and \be
\lim_{k\to\infty}Q(u_k)=\lim_{k\to\infty}\int_\Omega (|\nabla
u_k|^p+V|u_k|^p)\dx=0.\ee
 We say that a positive function $v\in
C^1_{\mathrm{loc}}(\Omega)$ is a {\em ground state} of the
functional $Q$ in $\Omega$ if $v$ is an
$L^p_{\mathrm{loc}}(\Omega)$ limit of a null sequence. If $Q\ge
0$, and $Q$ admits a ground state in $\Omega$, we say that $Q$ is
{\em degenerately positive} in $\Omega$.}
\end{definition}
\begin{definition} {\em The functional $Q$ is {\em nonpositive} in
$\Omega$ if it takes negative values on $\core$.}
\end{definition}
The main result of the present paper reads as follows:
\begin{theorem}
\label{main} Let $\Omega\subseteq\R^d$ be a domain, $V\in
L_\mathrm{loc}^\infty(\Omega)$, and $p\in(1,\infty)$. Suppose that
the functional $Q$ is nonnegative. Then
\begin{itemize}
\item[(a)] $Q$ has either a weighted spectral gap or a ground
state.

\item[(b)] If the functional $Q$ admits a ground state $v$, then $v$ satisfies
(\ref{groundstate}).

\item[(c)]
The functional $Q$ admits a ground state if and only if
(\ref{groundstate}) admits a unique positive supersolution.

\item[(d)]
 If $Q$ has a ground state $v$, then there exists a positive continuous function $W$ in
$\Omega$, such that for every $\psi\in C_0^\infty(\Omega)$
satisfying $\int \psi v \,\mathrm{d}x \neq 0$ there exists a
constant $C>0$ such that the following inequality holds:
 \be\label{Poinc}
 C^{-1}\int_{\Omega} W|u|^p\,\mathrm{d}x\le Q(u)+C\left|\int_{\Omega} \psi
 u\,\mathrm{d}x\right|^p\qquad
 \forall u\in C_0^\infty(\Omega).\ee
\end{itemize}
\end{theorem}
Theorem \ref{main} extends \cite[Theorem~1.5]{PT2} that deals with
the linear case $p=2$.

The outline of the present paper is as follows. The next section
gives preliminary results concerning positive solutions of the
equation $Q'(u)=0$. In particular, it introduces the generalized
Picone's identity \cite{AH1,AH2} which is a crucial tool in our
study. Section~\ref{sec:main} is devoted to the proof of
Theorem~\ref{main}. In Section~\ref{sec:crit}, we study, using
Theorem~\ref{main}, criticality properties of the functional $Q$
along the lines of criticality theory for second-order linear
elliptic operators \cite{P0,P2}.

In Section~\ref{sec:mingr} we prove, for $1<p\leq d$, the
existence of a unique (up to a multiplicative constant) positive
solution of the equation $Q'(u)=0$ in $\Omega\setminus \{x_0\}$
which has a minimal growth in a neighborhood of infinity in
$\Omega$. The proof of the above result relies on an unpublished
lemma of L.~V\'{e}ron (Lemma~\ref{lemveron}) concerning the exact
asymptotic behavior of a singular positive solution of the
equation $Q'(u)=0$ in a punctured neighborhood of $x_0$. We thank
Professor V\'{e}ron for kindly supplying the proof of this key
result. The study of positive solutions which has a minimal growth
in a neighborhood of infinity in $\Omega$ leads us to another
characterization of strict positivity in terms of these solutions.

In Section~\ref{sec:open}, we pose a number of open problems
suggested by the results of the present paper. The Appendix
contains a new energy estimate for the functional $Q$ valid for
$p> 2$.  This estimate leads to an alternative proof of
Lemma~\ref{someballs} for the case $p> 2$.

\mysection{Positive solutions and Picone identity}
Let $v>0$, $v\in C^1_{\mathrm{loc}}(\Omega)$, and $u\geq 0$,
$u\in\core$. Denote \be R(u,v):=|\nabla u|^p- \nabla
\left(\frac{u^{p}}{v^{p-1}}\right)\cdot|\nabla v|^{p-2}\nabla v,
\ee and
 \be L(u,v):=|\nabla u|^p+(p-1)\frac{u^p}{v^p}|\nabla
v|^p-p\frac{u^{p-1}}{v^{p-1}}\nabla u\cdot|\nabla v|^{p-2}\nabla
v.\ee Then the following {\em (generalized) Picone identity} holds
\cite{AH1,AH2} \be\label{picone0} R(u,v)=L(u,v). \ee
 Write
$L(u,v)=L_1(u,v)+L_2(u,v)$, where
 \be \label{L1}
L_1(u,v):=|\nabla u|^p+(p-1)\frac{u^p}{v^p}|\nabla
v|^p-p\frac{u^{p-1}}{v^{p-1}}|\nabla u||\nabla v|^{p-1},\ee
 and
 \be\label{L2}
 L_2(u,v):=p\frac{u^{p-1}}{v^{p-1}}|\nabla
v|^{p-2}(|\nabla u||\nabla v|-\nabla u\cdot\nabla v)\ge 0.\ee
 From the obvious inequality $t^p+(p-1)-pt\ge 0$, we also have that
$L_1(u,v)\ge 0$. Therefore, $L(u,v)\ge 0$ in $\Omega$.

Let $v\in C_{\mathrm{loc}}^1(\Omega)$ be a positive solution
(resp. supersolution) of (\ref{groundstate}). Using
(\ref{picone0}) and (\ref{solution}) (resp.
(\ref{supersolution})),
 we infer that for every $u\in\core$, $u\ge 0$,
 \be\label{QL} Q(u)=\int_\Omega L(u,v)\dx\geq 0, \qquad
 \left(\mbox{resp. } Q(u)\geq \int_\Omega L(u,v)\dx\geq 0\right).
 \ee

For any smooth subdomain $\Omega'\Subset\Omega$ consider the
variational problem \be \label{mu} \lambda_{1,p}(\Omega'):=\inf_{
u\in W_0^{1,p}(\Omega')}\dfrac{\int_{\Omega'}(|\nabla
u|^p+V|u|^p)\dx}{\int_{\Omega'} {|u|^p}\dx}\,.\ee
 It is well-known that for such a subdomain, (\ref{mu}) admits (up to a multiplicative
constant) a unique minimizer $\varphi$ \cite{DKN,GS}. Moreover,
$\varphi$ is a positive solution of the quasilinear eigenvalue
problem
\begin{equation}
  \begin{cases}
Q'(\varphi)=\lambda_{1,p}(\Omega')|\varphi|^{p-2}\varphi    & \text{ in } \Omega', \\
    \varphi=0 & \text{ on } \partial \Omega'.
  \end{cases}
\end{equation}
$\lambda_{1,p}(\Omega')$ and $\varphi$ are called the {\em
principal eigenvalue and eigenfunction of the operator $Q'$},
respectively.

\vskip 3mm

The following theorem was proved by J.~Garc\'{\i}a-Meli\'{a}n, and
J.~Sabina de Lis \cite{GS} (see also \cite{AH1,AH2}).
\begin{theorem}\label{thmGS} Assume that $\Omega\subset\R^d$ is a bounded $C^{1+\alpha}$-domain,
$0<\alpha<1$, and suppose that $V\in L^\infty(\Omega)$. Then the
following assertions are equivalent

\begin{itemize}
 \item[(i)] $Q'$  satisfies the
maximum principle: If $u$ is a solution of the equation
$Q'(u)=f\geq 0$ in $\Omega$ with some $f\in L^\infty(\Omega)$, and
satisfies $u\geq 0$ on $\partial\Omega$, then $u$ is nonnegative
in $\Omega$.

\item[(ii)] $Q'$  satisfies the strong maximum principle: If $u$
is a solution of the equation $Q'(u)=f\gneqq 0$ in $\Omega$ with
some $f\in L^\infty(\Omega)$, and satisfies $u\geq 0$ on
$\partial\Omega$, then $u>0$ in $\Omega$.

\item[(iii)] $\lambda_{1,p}(\Omega)>0$.

\item[(iv)] For some $0\lneqq f\in L^\infty(\Omega)$ there exists a positive strict
supersolution $v$ satisfying  $Q'(v)=f$ in $\Omega$, and $v=0$ on
$\partial\Omega$.

\item[(iv')] There exists a positive strict supersolution $v$
satisfying   $Q'(v)=f\gneqq 0$ in $\Omega$, such that $v\in
C^{1+\alpha}(\partial \Omega)$ and $f\in L^\infty(\Omega)$.

\item[(v)] For each nonnegative $f\in C^\alpha(\Omega)\cap
L^\infty(\Omega)$ there exists a unique weak nonnegative solution
of the problem $Q'(u)=f$ in $\Omega$, and $u=0$ on
$\partial\Omega$.
 \end{itemize}
 \end{theorem}
We shall need also the following comparison principle of
J.~Garc\'{\i}a-Meli\'{a}n, and J.~Sabina de Lis \cite[Theorem
5]{GS}.
\begin{theorem}\label{CP}
Let $\Omega\subset\mathbb{R}^d$ be a bounded domain of class
$C^{1,\alpha}$, $0< \alpha\leq 1$. Assume that
$\lambda_{1,p}(\Omega)>0$ and let $u_i\in W^{1,p}(\Omega)\cap
L^\infty(\Omega)$ satisfying $Q'(u_i)\in L^\infty(\Omega)$,
$u_i|_{\partial \Omega}\in C^{1+\alpha}(\partial \Omega)$, where
$i=1, 2$. Suppose further that the following inequalities are
satisfied
\begin{equation}
\left\{
\begin{array}{rclc}
 Q'(u_1)\!\!\!&\leq &\!\!\! Q'(u_2) \qquad &\mbox{ in }
\Omega,\\
Q'(u_2)\!\!\!&\geq &\!\!\! 0 \qquad &\mbox{ in }
\Omega,\\
u_{1}\!\!\!&\leq&\!\!\!u_{2} \qquad &\mbox{ on }\partial\Omega,\\
u_{2}\!\!\!&\geq&\!\!\!0 \quad&\mbox{ on } \partial \Omega.\\
\end{array}\right.
\end{equation}
Then
$$ u_{1}\leq u_{2} \qquad\mbox{ in }\Omega.$$
\end{theorem}

The following theorem generalizes the well-known
Allegretto-Piepenbrink theorem (see \cite[Theorem 2.12]{Cycon} and
the references therein).
\begin{theorem}\label{pos} Let Q be a functional of the form (\ref{Q}).
Then the following assertions are equivalent
\begin{itemize}
\item[(i)] The functional $Q$ is nonnegative on $C_0^\infty(\Omega)$.

\item[(ii)] Equation (\ref{groundstate}) admits a global positive
solution. \item[(iii)] Equation (\ref{groundstate}) admits a
global positive supersolution.
\end{itemize}
\end{theorem}
\begin{proof}
(i) $\Rightarrow$ (ii): Assume that $Q\geq 0$ on
$C_0^\infty(\Omega)$. Then $Q$ is nonnegative on
$C_0^\infty(\Omega')$ for any smooth bounded domain
$\Omega'\Subset\Omega$. Fix an {\em exhaustion}
$\{\Omega_{N}\}_{N=1}^{\infty}$ of $\Omega$ (i.e., a sequence of
smooth, relatively compact domains such that $x_0\in \Omega_1$,
$\mbox{cl}({\Omega}_{N})\subset \Omega_{N+1}$ and
$\cup_{N=1}^{\infty}\Omega_{N}=\Omega$). By the strict
monotonicity of $\lambda_{1,p}(\Omega)$ as a function of $\Omega$
\cite[Theorem 2.3]{AH1}, it follows that
$\lambda_{1,p}(\Omega_N)>0$ for all $N\geq 1$.  Let $f_N\in
C_0^\infty(\Omega_N\setminus\Omega_{N-1})$ be a nonnegative
nonzero function. By Theorem~\ref{thmGS}, there exists a unique
positive solution of the problem
$$Q'(u_N)=f_N \; \mbox{ in }
\Omega_N,  \qquad u_N=0 \;\mbox{ on }\partial\Omega_N.$$ Set
$v_N(x):=u_N(x)/u_N(x_0)$. By Harnack's inequality and elliptic
regularity \cite{Serrin1,T}, it follows that $\{v_N\}$ admits a
subsequence which converges locally uniformly to a positive
solution $v$ of the equation $Q'(u)=0$ in $\Omega$.

(ii) $\Rightarrow$ (iii) is trivial.

(iii) $\Rightarrow$ (i): Suppose that (\ref{groundstate}) admits a
global positive supersolution $v$ in $\Omega$. If $u\in
C_0^\infty(\Omega)$ is a nonnegative function, then by (\ref{QL})
$Q(u)\geq 0$. Let $u\in\sobloc$ be a nonnegative function with
compact support. By taking a sequence $\{u_k\}\subset\core$ of
nonnegative functions such that $u_k\to u$ in $\sobloc$, we infer
that $Q(u)\geq 0$. Thus, $Q\geq 0$ on the cone of all $\sobloc$
nonnegative functions with compact support in $\Omega$. Since
$Q(u)=Q(|u|)$ on $\core$, it follows that $Q$ is nonnegative on
$\core$.
\end{proof}

\mysection{Proof of Theorem~\ref{main}} \label{sec:main}
We shall start with the following two lemmas. For $B\subset\Omega$
let

\be \label{cb} c_B:=\inf_{\stackrel{u\in\core}{\int_B |u|^p\dx=1}}
Q(u)=\inf_{\stackrel{0\leq u\in\core}{\int_B |u|^p\dx=1}} Q(u)
.\ee

\begin{lemma}
\label{allballs} If for every open set $B\Subset \Omega$, $c_B>0$,
then there exists a $W\in C(\Omega)$, $W>0$, such that \be
\label{sg} Q(u)\ge \int_\Omega W(x)|u(x)|^p\dx\qquad \forall
u\in\core.\ee In other words, $Q$ has a weighted spectral gap with
the weight $W$.
\end{lemma}
\begin{proof} Let $\{B_{j}\}_{j=1}^\infty$ be a
locally finite covering of $\Omega$ by open balls
$B_{i}\Subset\Omega$, and let $\{\chi_j\}$ be a locally finite
partition of unity on $\Omega$ subordinated to this covering. Set
$C_j:=\min\{c_{B_{j}},1\}$. Then
 \be\label{2j} 2^{-j} Q(u)\ge
2^{-j}c_{B_{j}} {\int_{B_{j}} |u|^p \dx}\geq 2^{-j}C_{j}
{\int_\Omega \chi_j|u|^p \dx}\qquad \forall u\in\core. \ee
By summation (\ref{2j}) over
$j\in\N$ and by interchanging the order of summation and integration
we obtain (\ref{sg}) with
 $W(x):=\sum_{j=1}^\infty 2^{-j}C_j\chi_j(x)$.
\end{proof}

\begin{lemma}
\label{someballs} If there exists a nonempty open set
$B\Subset\Omega$ such that $c_B=0$, then $Q$ admits a ground
state.
\end{lemma}

\begin{proof} Fix a positive (super)solution $v\in
C^1(\Omega)$. Since $c_B=0$, there exists a sequence
$\{u_k\}\subset\core$, $u_k\ge 0$, such that $\int_B|u_k|^p\dx=1$
and $Q(u_k)\to 0$.  Let $\omega\Subset \Omega$ be an open
connected set containing $B$.

\vskip 2mm

\noindent{\bf Step 1.} Let $\omega^\prime\subset\omega$. By
(\ref{QL}), and since $L(u_k,v)\ge 0$, we have that
 \bea
\int_{\omega^\prime}|\nabla u_k|^p\dx\le \int_\Omega L(u_k,v)\dx+
\int_{\omega^\prime}C u_k^{p-1}|\nabla u_k|\dx
\le \nonumber\\
o(1)+C\int_{\omega^\prime}u_k^{p-1}|\nabla u_k|\dx, \eea with the
constant $C=C(\omega)$ independent of $\omega^\prime$.
Invoking Young's inequality, we arrive at
\bea \int_{\omega^\prime}|\nabla u_k|^p\dx\le
o(1)+\frac12\int_{\omega^\prime}|\nabla
u_k|^p\dx+C\int_{\omega^\prime} u_k^p\dx. \eea
Therefore,
\be\label{eq1} \int_{\omega^\prime}|\nabla u_k|^p\dx\le
\\C(\omega)\int_{\omega^\prime}u_k^p\dx +o(1). \ee

\vskip 2mm

\noindent{\bf Step 2.} Let
\be\omega_0:=\{x\in\omega:\;\exists\rho(x)\in(0,
d(x,\Omega\setminus\omega)),\;
\sup_k\int_{B_{\rho(x)}(x)}|u_k|^p\dx<\infty\}.\ee
Since $\int_B|u_k|^p\dx=1<\infty$, $B\subset\omega_0$. Moreover,
$\omega_0$ is an open set. Indeed, let $x\in\omega_0$,
$B_{\rho(x)}(x)\subset\omega_0$. Then for every point $y\in
B_{\rho(x)}(x)$ there is a $\rho(y)>0$ such that
$B_{\rho(y)}(y)\subset B_{\rho(x)}(x)$. Therefore,
$\int_{B_{\rho(y)}(y)}|u_k|^p\dx\le\int_{B_{\rho(x)}(x)}|u_k|^p\dx$.
Consequently $y\in \omega_0$, and $\omega_0$ is open.

Let us show now that $\omega_0$ is a relatively closed set in
$\omega$. We shall use the following version of Poincar\'e
inequality
 \be\label{poinver}
\int_{B_1(0)}\!\!\!|u|^p\dx\le C\!\! \int_{B_1(0)}\!\!|\nabla
u|^p\dx+C_r \left|\int_{B_r(0)}\!\!\!u\dx\right|^p \quad \forall
r\in(0,1),\; \forall u\in C^\infty(\R^d),\ee
 which follows for example
from \cite[Theorem 4.2.1]{Z}.
 It easily follows from (\ref{poinver}) that for
every $\epsilon>0$ and $\rho>0$ there exist $\delta_{\epsilon}>0$
and $C(\epsilon,\rho)$ such that for every $x\in\omega$,
$\delta\in(0,\min\{\delta_\epsilon,d(x,\partial\omega)\})$, and
$u\in C_0^\infty(\Omega)$,
\be \label{scaledP} \int_{B_\delta(x)}|u|^p\dx\le \epsilon
\int_{B_\delta(x)}|\nabla
u|^p\dx+C(\epsilon,\rho)\int_{B_\rho(x)}|u|^p\dx.\ee

Let $x_j\in\omega_0$, $x_j\to x_0\in\omega$. Let $\epsilon<
(C(\omega))^{-1}/2$, where $C(\omega)$ is the constant in
(\ref{eq1}). Let $\delta_\epsilon>0$ be as in (\ref{scaledP}) and
fix $\delta\in(0,\min\{\delta_\epsilon,d(x,\partial\omega)\})$.
Finally, choose $j$ such that $|x_0-x_j|\le \frac{\delta}{2}$.
Then, with $\rho=\rho(x_j)$, (\ref{eq1}) and (\ref{scaledP}) imply
\be \label{eq310}\frac12C(\omega)^{-1}\int_{B_\delta(x_j)}|\nabla
u_k|^p\dx\le
C(\epsilon,\rho(x_j))\int_{B_{\rho(x_j)}(x_j)}|u_k|^p\dx+o(1).\ee
The right hand side of (\ref{eq310}) is bounded in $k$ by the
definition of $\omega_0$. Thus $\int_{B_\delta(x_j)}|\nabla
u_k|^p\dx$ is bounded, and by (\ref{scaledP}),
$\int_{B_\delta(x_j)}|u_k|^p\dx$ is also bounded. Since, by the
choice of $j$, $B_{\delta/2}(x_0)\subset B_\delta(x_j)$, it
follows that $\int_{B_{\delta/2}(x_0)}(|\nabla u_k|^p+|u_k|^p)\dx$
is bounded and consequently, $x_0\in\omega_0$, and $\omega_0$ is
also relatively closed.

Since $\omega$ is connected and $\omega_0\neq\emptyset$, it
follows that $\omega_0=\omega$ and $\{u_k\}$ is bounded in
$L^{p}_\mathrm{loc}(\Omega)$. Invoking again (\ref{eq1}) it
follows that $\{u_k\}$ is bounded in
$W^{1,p}_\mathrm{loc}(\Omega)$.

\vskip 2mm

\noindent{\bf Step 3.} Consider now a weakly convergent renamed
subsequence $u_k\rightharpoonup u$ in
$W^{1,p}_\mathrm{loc}(\Omega)$. Let $\omega\Subset \Omega$ be a
smooth domain, and set $Q^\omega(u):=\int_\omega L(u,v)\dx$. We
claim that the functional $Q^\omega(u)$ is weakly lower
semicontinuous in $W^{1,p}(\omega)$. Indeed, the functionals
$\int_\omega |\nabla u|^p\dx$ and $\int_\omega
(p-1)\frac{|u|^p}{v^p}|\nabla v|^p\dx$ are weakly lower
semicontinuous in $W^{1,p}(\omega)$ since their Lagrangians
($\mathcal{L}(q,z,x)=|q|^p$ and
$\mathcal{L}(q,z,x)=(p-1)\frac{|z|^p}{v(x)^p}|\nabla v(x)|^p$,
resp.)  are convex functions of $q$. So, it suffices to show that
the functional \be J^\omega(u):=\int_\omega
\frac{u^{p-1}}{v^{p-1}}\nabla u\cdot|\nabla v|^{p-2}\nabla v\dx\ee
is weakly continuous on any sequence $\{u_k\}$ satisfying
$u_k\rightharpoonup u$ in $W^{1,p}(\omega)$. Indeed,
\bea J^\omega(u_k)-J^\omega(u)=\nonumber\int_\omega
|\nabla v|^{p-2}v^{1-p}\nabla v\cdot \nabla u_k(u_k^{p-1}-u^{p-1})\dx+\\
\int_\omega \nabla (u_k-u)\cdot\nabla v |\nabla
v|^{p-2}v^{1-p}u^{p-1}\dx.\label{J}\eea
Consider the first term of the right hand side of (\ref{J}). Since
$u_k\rightharpoonup u$ in $W^{1,p}(\omega)$, it follows by the
compactness of the local Sobolev imbeddings that (up to a
subsequence)  $u_k\to u$ in $L^{p}(\omega)$. Then, for a renamed
subsequence, there exists a $U\in L^{p}(\omega)$, such that $0\le
u_k\le U$ and $u_k\to u$ a.e. in $\omega$. Therefore,
$u_k^{p-1}\le U^{p-1}\in L^{p'}(\omega)$, where $p':=p/(p-1)$ is
the conjugate exponent of $p$. Consequently,
 \be\label{eq14}|u_k^{p-1}-u^{p-1}|^{p'}\leq C(U^p+u^p)\in L^1(\omega).
  \ee
 Hence by H\"older's inequality and Lebesgue's dominated
  convergence theorem, \be \left|\int_\omega
|v|^{p-2}v^{1-p}\nabla v \cdot \nabla u_k (u_k^{p-1}-u^{p-1})
\dx\right|\le C\|\nabla u_k\|_p\|u_k^{p-1}-u^{p-1}\|_{p'}\to 0.\ee

Consider the functional $$\Phi (w):=\int_\omega \nabla
w\cdot\nabla v |\nabla v|^{p-2}v^{1-p}u^{p-1}\dx.$$ Note $|\nabla
v|^{p-1}v^{1-p}u^{p-1}\in L^{p'}(\omega)$. Therefore, H\"older's
inequality implies that $\Phi$ is a continuous functional on
$W^{1,p}(\omega)$. Hence, by the definition of weak convergence,
the second term of the right hand side of (\ref{J}) converges to
zero.

We conclude that $0\le Q^\omega(u)\le\liminf Q^\omega(u_k)=0$.
Moreover, $\int_B u^p\dx=1$. Now we repeat the argument of
\cite{AH1}. Since $Q^\omega(u)=0$ for every subdomain
$\omega\Subset \Omega$ containing $B$, it follows that
$L_1(u,v)=0$ and $L_2(u,v)=0$. Recall that $f(t)=t^p+p-1-pt$ is a
nonnegative function on $\mathbb{R}_+$ which attains its zero
minimum only at $t=1$. Therefore, $L_1(u,v)=0$ implies that
$u^{-1}|\nabla u|=v^{-1}|\nabla v|$. On the other hand,
$L_2(u,v)=0$ implies that $\nabla u$ is parallel to $\nabla v$.
Hence, $u=cv$, where $c>0$. The value of $c$ is determined by the
condition $\int_B (cv)^p\dx =1$. Therefore $cv$ is the limit of
every weakly convergent subsequence of $\{u_k\}$. It follows that
the original sequence $u_k\to cv$ in $L^p_\mathrm{loc}(\Omega)$.
\end{proof}

\begin{remark}
\label{wlsc}{\em
 The argument in Step 3 of the proof of Lemma
\ref{someballs} shows that the functional $Q$ has some weakly
lower semicontinuity properties, as Proposition~\ref{prwlsc} below
demonstrates.
 }\end{remark}
First, we need to extend the definition of the
functional $Q$ from $\core$ to a larger set in $\sobloc$.
\begin{definition}{\em
 Let $v$ be a positive solution of the equation $Q'(u)=0$ in
$\Omega$. We define a functional $Q_v:\sobloc\to[0,+\infty]$  by
$$Q_v(u):=
  \begin{cases}
   \int_\Omega L(|u|,v)\dx & \text{if this integral is finite}, \\
    \infty & \text{otherwise}.
  \end{cases}
  $$
 Also, for $u\in \sobloc$ and any
$\omega\Subset\Omega$ define $Q^\omega_v(u):=\int_\omega
L(|u|,v)\dx$. }\end{definition}
We have

\begin{proposition}\label{prwlsc} Let
\bea\nonumber
 D:=\left\{ u\in\sobloc\mid Q_v(u)<\infty,
\mbox{ and } \exists\{u_k\}\subset
\cup_{\omega\Subset\Omega}W^{1,p}_0(\omega)
\mbox{ s.t.}\quad\right.\\
 \left. u_k\to u \mbox{ in } \sobloc, \mbox{ and } Q_v(u_k)\to
Q_v(u)\right\}.\label{defD}
 \eea
 Then the functional $Q$ on $\core$ admits an extension to the set $D$
given by $Q(u):=Q_v(u)$. This extension is independent of the
positive solution $v$.
Moreover, the functional $Q$ is continuous on
$W^{1,p}_0(\omega)\subset D$ for every $\omega\Subset\Omega$, and
is weakly lower semicontinuous in the following sense: \be
\label{weakcont} u_k,u\in D,\; u_k\rightharpoonup u \text{ in }
\sobloc,\;  \sup_{k\in\mathbb{N}} Q(u_k)<\infty \Rightarrow
\liminf_{k\to\infty} Q(u_k)\ge Q(u). \ee
\end{proposition}

Note that the functional $Q^{1/p}$ is generally not a norm (see
the discussion prior to Problem~\ref{pr:conv}). If $Q^{1/p}$ is a
norm and $Q$ is strictly positive, then the set $D$ is a closure
of $\core$ in that norm and so it is a Banach space. In this case,
the condition in (\ref{weakcont}) is equivalent to weak
convergence in $D$, and $Q$ is weakly lower semicontinuous as a
monotone increasing function of the norm.

\begin{proof}
The functionals $Q$ and $Q_v$ obviously admit continuous
extensions to  $W^{1,p}_0(\omega)$ for every
$\omega\Subset\Omega$. Since $Q$ and $Q_v$ are even  and  coincide
on nonnegative $\core$-functions, their respective extensions to
$W^{1,p}_0(\omega)$ are equal. Consequently, the set  $D$  is
independent of $v$. Moreover, the functional $Q_v$  evaluated on
$D$ is independent of $v$ and thus defines the extension of $Q$.

Let $u_k\rightharpoonup u$ in $\sobloc$, $u_k,u\in D$, and
$Q_v(u_k)\le C$ for some $C>0$. Assume first that $u_k\ge 0$. Then
$Q_v(u_k)\ge Q^\omega_v(u_k)$  for every open
$\omega\Subset\Omega$. Step 3 of Lemma~\ref{someballs} implies
that $\liminf Q_v(u_k)\ge\liminf Q^\omega_v(u_k)\ge
Q^\omega_v(u)$. Since $\omega$ is arbitrary, we have $\liminf
Q_v(u_k)\ge Q_v(u)=Q(u)$.

Let now remove the restriction $u_k\ge 0$. Let
$\omega\Subset\Omega$ be a smooth domain. Due to the compactness
of the imbedding of $W^{1,p}(\omega)$ into $L^p(\omega)$,  $u_k\to
u$ a.e. on a renamed subsequence,
 and
$|u_k|\le U$ with some $U\in L^p_{\mathrm{loc}}(\Omega)$.
Consequently $|u_k|\to |u|$ in $L^p_{\mathrm{loc}}(\Omega)$. Since
$|u_k|$ is bounded with respect to the seminorms of $\sobloc$,
$|u_k|\rightharpoonup |u|$ in $\sobloc$. Since the functional
$Q_v$ is even, by the previous argument it follows that
$$\liminf_{k\to\infty}
Q_v(u_k)=\liminf_{k\to\infty} Q_v(|u_k|)\ge
Q_v(|u|)=Q_v(u)=Q(u),$$ and the proposition is proved.
\end{proof}

{\it Proof of Theorem~\ref{main}.} Part (a) follows from
Lemma~\ref{allballs} and Lemma~\ref{someballs}.

\noindent To prove (b) and (c), observe that from
Lemma~\ref{someballs} it follows that for each positive
supersolution $v$ of (\ref{groundstate}), any null sequence
$\{u_k\}$ converges to a constant multiple of $v$. This implies
that all positive supersolutions of (\ref{groundstate}) are scalar
multiples of the same function. On the other hand, if $Q$ admits a
weighted spectral gap with a weight $W$, then by
Theorem~\ref{pos}, the equation $\frac{1}{p}Q'(u)-W|u|^{p-2}u=0$
admits a positive solution $v$ in $\Omega$. So, $v$ is a strictly
positive supersolution of the equation $Q'(u)=0$ in $\Omega$. In
addition, by Theorem~\ref{pos}, the equation $Q'(u)=0$ admits a
positive solution $w$ in $\Omega$. Clearly $v\neq w$, and the
equation $Q'(u)=0$ admits two linearly independent positive
supersolutions in $\Omega$.

It remains to prove (d). First we claim that for every open set
$B\Subset\Omega$ there is a strictly positive continuous function
$W$ such that

\be \label{2.5} \int_\Omega W|u|^p\dx\le Q(u)+\int_B|u|^p\dx\qquad
\forall u\in\core .
 \ee

Indeed, denote the functional in the right hand side of
(\ref{2.5}) by $\tilde{Q}$. Clearly, $\tilde{Q}\geq 0$ on $\core$.
Suppose that $\tilde{Q}$ admits a null sequence $\{u_k\}$, then
$\{u_k\}$ is a null sequence of $Q$, which implies
$\int_B|u_k|^p\dx\to \int_B|v|^p\dx$, where $v$ is a ground
state for $Q$. Consequently, $\liminf \tilde{Q}(u_k)>0$, which
contradicts the definition of $\{u_k\}$. Therefore, part (a) of
the present theorem implies that $\tilde{Q}$ admits a weighted
spectral gap.

Therefore, in order to prove (\ref{Poinc}), it suffices to show
that for some open $B\Subset\Omega$,

\be \label{Poinc2} \int_B|u|^p\dx\le C\left(Q(u)+\left|\int_\Omega
u\psi \dx\right|^p\right) \qquad \forall u\in\core. \ee

Suppose that (\ref{Poinc2}) fails. Then there is a sequence
$\{u_k\}\subset \core$ such that $\int_B|u_k|^p\dx=1$, $Q(u_k)\to
0$ and $\int_\Omega u_k\psi \dx\to 0$. Since $\{u_k\}$ is a null
sequence it converges in $L^p_{\mathrm{loc}}(\Omega)$ to $v$,
where $v>0$ is the ground state of $Q$. Then $\int_\Omega u_k\psi
\dx\to \int_\Omega v\psi \dx\neq 0$, and we arrive at a
contradiction.
 \qed

 For $u\in \core$, we define \bean\tilde{Q}(u):=
  \begin{cases}
    Q(u) & \text{ if Q has a weighted spectral gap,} \\
    Q(u)+C\left|\int_{\Omega} \psi
 u\,\mathrm{d}x\right|^p & \text{ if Q has a ground state,}
  \end{cases}
 \eean
where $C$ is the constant in (\ref{Poinc}).
\begin{proposition}
\label{compact}  For any $C>0$, the set $S:=\{u\in \core \mid
u\geq 0,\; \tilde{Q}(u)\le C \}$ is bounded in $\sobloc$ and
therefore, it is relatively compact in
 $L^p_\mathrm{loc}(\Omega)$.
\end{proposition}

\begin{proof} Let $u\in S$ and let $v$ be a positive solution of
(\ref{groundstate}). Let $\omega\Subset \Omega$ be an open set.
Then by (\ref{QL}), we have that \bean \int_\omega(|\nabla u|^p
+|u|^p)\dx\le\int_\omega (L(u,v)+C |u|^{p-1}|\nabla u|+|u|^p)\dx
\le \\
C+C\int_\omega |u|^{p-1}|\nabla u|\dx+\int_\omega |u|^p\dx. \eean
 In light of Young's inequality we obtain \bean \int_\omega(|\nabla u|^p
+|u|^p)\dx\le C+\frac12\int_\omega|\nabla u|^p\dx+C\int_\omega
|u|^p\dx. \eean
 Consequently,
\be\label{eq7} \int_\omega(|\nabla u|^p +|u|^p)\dx\le
\\C+C\int_\omega |u|^p\dx. \ee
 By (\ref{wsg}) or (\ref{Poinc}) and the definition of $S$, the right hand side of (\ref{eq7})
is uniformly bounded in $S$, and therefore $S$ is a bounded set in
$\sobloc$.
\end{proof}

\mysection{Criticality theory}\label{sec:crit}

In this section we prove several positivity properties of the
functional $Q$ along the lines of criticality theory for
second-order linear elliptic operators \cite{P0,P2}.

Recall that $Q$ is said to be  strictly positive in $\Omega$ if
$Q$ has a weighted spectral gap on $\core$. The functional $Q$ is
degenerately positive in $\Omega$ if $Q\ge 0$ on $\core$ and $Q$
admits a ground state in $\Omega$. The functional $Q$ is
nonpositive in $\Omega$ if it takes negative values on $\core$.
For $V\in L^\infty_{\mathrm{loc}}(\Omega)$, we denote in the
present section \be Q_V(u):=\int_\Omega(|\nabla u|^p+V|u|^p)\dx\ee
to emphasize the dependence of $Q$ on the potential $V$.
\begin{proposition}
\label{monPot} Let  $V_i\in L^\infty_{\mathrm{loc}}(\Omega)$ and
suppose that $V_2\gneqq V_1$. If  $Q_{V_1}\ge 0$ on $\core$, then
$Q_{V_2}$ is strictly positive in $\Omega$, and if $Q_{V_2}$ is
degenerately positive in $\Omega$, then $Q_{V_1}$ is nonpositive
in $\Omega$.
\end{proposition}
\begin{proof} Obviously,
$$Q_{V_2}(u)=Q_{V_1}(u)+\int(V_2-V_1)|u|^p\dx\ge 0\qquad \forall u\in \core.$$ Suppose that
$Q_{V_2}$ has a null sequence $\{u_k\}$ with a ground state $v$,
such that $u_k\to v$ in $L^p_{\mathrm{loc}}(\Omega)$.  Evaluating
the limit of $Q_{V_2}(u_k)$, we have by Fatou's lemma,\be
0=\lim_{k\to\infty} Q_{V_2}(u_k)\ge\liminf_{k\to\infty}
\int_\Omega(V_2-V_1)|u_k|^p \dx\ge \int_\Omega(V_2-V_1)|v|^p\dx>0,
\ee and we arrive at a contradiction.
\end{proof}

\begin{proposition}
\label{monDom}
Let $\Omega_1\subset\Omega_2$ be domains in $\R^d$ such that $\Omega_2\setminus\overline{\Omega_1}\neq\emptyset$.
Let $Q_V$ be defined on $C_0^\infty(\Omega_2)$.

1. If $Q_V\ge 0$ on $C_0^\infty(\Omega_2)$, then $Q_V$ is strictly
positive in $\Omega_1$.

2. If $Q_V$ is degenerately positive in $\Omega_1$, then $Q_V$ is
nonpositive in $\Omega_2$.
\end{proposition}
\begin{proof}
1. If $Q_V$ is strictly positive in $\Omega_2$, then the first
assertion is trivial. Suppose that $Q_V$ is degenerately positive
in $\Omega_2$, and let $v$ be the ground state of $Q_V$ in
$\Omega_2$. Take $\psi\in
C_0^\infty(\Omega_2\setminus\overline{\Omega_1})$ such that
$\int_{\Omega_2} v\psi\dx\neq 0$.

Due to (\ref{Poinc}) restricted to $u\in C_0^\infty(\Omega_1)$, we
conclude that $Q_V$ is strictly positive in $\Omega_1$.

2. Assume that $Q_V\ge 0$ on $C_0^\infty(\Omega_2)$. Then by the
first part, $Q_V$ is strictly positive in $\Omega_1$, which is a
contradiction.
\end{proof}

\begin{proposition}\label{Prop2}
 Let $V_0, V_1\in L^\infty_{\mathrm{loc}}(\Omega)$,
$V_0\neq V_1$. For $t\in \mathbb{R}$ we denote \be
Q_t(u):=tQ_{V_1}(u)+(1-t)Q_{V_0}(u),\ee and suppose that
$Q_{V_i}\geq 0$ on $\core$ for $i=0,1$.

Then the functional $Q_t\geq 0$ on $\core$ for all $t\in[0,1]$.
Moreover, if $V_0\neq V_1$, then $Q_t$ is strictly positive in
$\Omega$ for all $t\in(0,1)$.
\end{proposition}
\begin{proof} The first assertion is immediate. To prove the
second assertion, assume first that at least one of
$Q_{V_0},Q_{V_1}$, say, $Q_{V_0}$, is strictly positive with a
weight $W_0$. Then for $0\leq \tau<1$, the functional $Q_\tau$ is
strictly positive with the weight $(1-\tau) W_0$.

Assume now that both $Q_{V_0}$ and $Q_{V_1}$ are degenerately
positive with ground states $v_0, v_1$, respectively, and assume
that for some $\tau\in(0,1)$, $Q_{\tau}$ has a null sequence
$\{u_k\}$ and a ground state $v_\tau$.

Note that $v_\tau$ is not a multiple of $v_0$ or of $v_1$ since
$V_0\neq V_1$. Then there exist $\psi_i\in C_0^\infty(\Omega)$,
$i=0,1$, such  that  \be \label{orth} \int_\Omega\psi_i v_i\dx
\neq 0, \quad \mbox{ and } \int_\Omega\psi_i v_\tau\dx=0\quad
i=0,1.\ee
 By (\ref{Poinc}), for $i=0,1$
there exist a continuous function $W_i>0$ in $\Omega$, and a
constant $C_i>0$, such that
\be \int_\Omega W_i|u|^p\dx\le Q_i(u)+C_i\left|\int_\Omega
u\psi_i\dx\right|^p \qquad \forall u\in \core. \ee Let
$W_{\tau}:=\tau W_1+(1-\tau)W_0$. Then
\be \int_\Omega\!\! W_{\tau}|u|^p\dx \le
Q_{\tau}(u)+C_1\tau\left|\int_\Omega\!
u\psi_1\dx\right|^p\!+C_0(1-\tau)\left|\int_\Omega\!
u\psi_0\dx\right|^p \quad \forall u\!\in\! \core. \ee
Substituting $u=u_k$ and passing to the limit, taking into account
that $u_k\to v_\tau$ in $L^p_\mathrm{loc}(\Omega)$ as well as
(\ref{orth}) and Fatou's lemma, we have
\be 0<\int_\Omega W_{\tau}|v_\tau|^p\dx \le
C_1\tau\left|\int_\Omega
v_\tau\psi_1\dx\right|^p+C_0(1-\tau)\left|\int_\Omega
v_\tau\psi_0\dx\right|^p=0, \ee
and we arrive at a contradiction. Therefore, $Q_\tau$ does not
admit a ground state, and by Theorem~\ref{main}, the functional
$Q_\tau$ is strictly positive in $\Omega$.
\end{proof}
\begin{proposition}\label{strictpos}
Let $Q_V$ be a strictly positive functional in $\Omega$. Consider
$V_0\in L^\infty(\Omega)$ such that $V_0\ngeq 0$ and $\supp
V_0\Subset\Omega$. Then there exist $\tau_+>0$ and $-\infty\leq
\tau_-<0$ such that $Q_{V+tV_0}$ is strictly positive in $\Omega$
for $t\in(\tau_-,\tau_+)$, and $Q_{V+\tau_+ V_0}$ is degenerately
positive in $\Omega$.
\end{proposition}
\begin{proof}
If $Q_V$ has a weighted spectral gap in $\Omega$ with a continuous
weight $W$, then $Q_{V+tV_0}$ satisfies the inequality \be
\int_\Omega (W+tV_0)|u|^p\dx \le Q_{V+tV_0}(u) \qquad \mbox{ on }
\core. \ee The weight $W+tV_0$ is strictly positive in $\Omega$
for $|t|$ small, since $W$ is a strictly positive continuous
function and $V_0$ is a bounded function with compact support. By
Proposition~\ref{Prop2}, the set of $t\in\R$, for which
$Q_{V+tV_0}$ is strictly positive, is an interval. Moreover, this
interval does not extend to $+\infty$. Indeed, take $u_0\in
C_0^\infty$ such that $\int_\Omega V_0|u_0|^p\dx<0$ to see that
$Q_{V+tV_0}(u_0)<0$ for $t$ sufficiently large. Note that by the
above argument, this interval extends to $-\infty$ if and only if
$V_0\leq 0$.

Let $\tau_+$ be the right endpoint of this interval. Obviously
$Q_{V+\tau_+ V_0}$ is nonnegative on $\core$. If, on the other
hand, the functional $Q_{V+\tau_+ V_0}$ is strictly positive in
$\Omega$, then by the preceding argument, there exists $\delta>0$
such that the functional $Q_{V+(\tau_+\, + \delta) V_0}$ is
strictly positive in $\Omega$, which contradicts the definition of
$\tau_+$.
\end{proof}
\begin{proposition}\label{propintcond}
Let $Q_V$ be a degenerately positive functional in $\Omega$, and
let $v$ be the corresponding ground state. Consider $V_0\in
L^\infty(\Omega)$ such that $\supp V_0\Subset\Omega$. Then there
exists $0<\tau_+\leq\infty$ such that $Q_{V+tV_0}$ is strictly
positive in $\Omega$ for $t\in(0,\tau_+)$ if and only if
 \be\label{intcond}
  \int_\Omega
V_0|v|^p\dx>0.
 \ee
\end{proposition}
\begin{proof}
Suppose that there exists $t>0$ such that $Q_{V+tV_0}$ is strictly
positive in $\Omega$. Then there exists $W\in C(\Omega)$, $W>0$,
such that
 \be \label{sg7} Q_V(u)+t\int_\Omega V_0|u|^p\dx\ge
\int_\Omega W|u|^p\dx\qquad \forall u\in\core.\ee
Let $\{u_k\}$ be a null sequence for the functional $Q_V$, and let
$v>0$ be the ground state  of $Q_V$ which is the
$L^p_{\mathrm{loc}}(\Omega)$ limit of $\{u_k\}$. By (\ref{sg7})
and Fatou's lemma we have \bean \label{sg8} t\int_\Omega
V_0|v|^p\dx=\lim_{k\to\infty}Q_V(u_k)+t\lim_{k\to\infty}\int_\Omega
V_0|u_k|^p\dx\ge\\ \liminf_{k\to\infty}\int_\Omega W|u_k|^p\dx\ge
\int_\Omega W|v|^p\dx>0.\eean Thus, \eqref{intcond} is satisfied.

 Suppose that (\ref{intcond}) holds true, but for any $t>0$ the
functional $Q_{V+tV_0}$ is nonpositive in $\Omega$. Therefore, for
any $t>0$ there exists $u_t\in \core$ such that
 \be\label{neg}
 Q_V(u_t)+t\int_\Omega V_0|u_t|^p\dx<0.
 \ee
Clearly, we may assume that $u_t\geq 0$. Since $Q_V(u_t)\ge 0$ it
follows that
 \be\label{neg5}
 \int_\Omega V_0|u_t|^p\dx<0.
 \ee
In particular, $\supp(u_t)\cap \supp(V_0)\neq \emptyset$.
Therefore, we may assume that
$$\int_{\supp(V_0)}
|u_t|^p\dx=1.$$ It follows that
 \be
 \lim_{t\to 0}t\int_\Omega V_0|u_t|^p\dx=0,
 \ee
and by the nonnegativity of $Q_V$ and (\ref{neg}),
 \be
 0\le \liminf _{t\to 0} Q_V(u_t)\le\limsup_{t\to 0} Q_V(u_t)\le 0.
 \ee
It follows that $\{u_t\}$ is a null sequence, and therefore,
$u_t\to v$ in $L^p_{\mathrm{loc}}(\Omega)$ as $t\to 0$, where $v$
is the corresponding ground state  of $Q_V$. Using a standard
argument similar to (\ref{eq14}) we have (for a subsequence)
 \be\label{eq12}\lim_{t\to 0}\int_\Omega
V_0|u_t|^p\dx=\int_\Omega V_0|v|^p\dx.
 \ee
  Combining (\ref{intcond}),
and (\ref{neg5}) and (\ref{eq12}), we obtain
 \be
0<\int_\Omega V_0|v|^p\dx=\lim_{t\to 0}\int_\Omega
V_0|u_t|^p\dx\leq 0,
  \ee
which is a contradiction.
\end{proof}
\begin{remark}
{\em An alternative proof of Proposition~\ref{propintcond} can be
derived from the (nonsymmetric) technique in \cite{pcr}.}
\end{remark}
\mysection{Minimal growth}\label{sec:mingr}

In this section we study the existence of positive solutions of
the equation $Q'(u)=0$ of minimal growth in a neighborhood of
infinity in $\Omega$, and obtain a new characterization of strict
positivity in terms of these solutions.

Throughout this section we assume that $1<p\leq d$. Therefore, for
any $x_0\in \Omega$, any positive solution $v$ of the equation
$Q'(u)=0$ in a punctured neighborhood of $x_0$ has either a
removable singularity at $x_0$, or
\begin{equation}\label{nonremov}
  v(x)\asymp\begin{cases}
    |x-x_0|^{\alpha(d,p)} & p<d, \\
     -\log |x-x_0| & p=d,
  \end{cases} \qquad \mbox{ as } x\to x_0,
\end{equation}
where $\alpha(d,p):=(p-d)/(p-1)$ \cite{Serrin1,Serrin2,V}. Here
$f\asymp g$ means that $c\leq f/g\leq C$, where $c$ and $C$ are
positive constants. In particular, in the nonremovable case,
\begin{equation}\label{nonremov1}
\lim_{x\to x_0} v(x)=\infty.
\end{equation}
\begin{lemma}[L.~V\'{e}ron, private communication]\label{lemveron}
Assume that $1<p\leq d$, and let $x_0\in \mathbb{R}^d$ be fixed.
Suppose that $v$ is a positive solution of the equation $Q'(u)=0$
in a punctured neighborhood of $x_0$ which has a nonremovable
singularity at $x_0$. Then
 \begin{equation}\label{nonremovasymp}
  v(x)\sim\begin{cases}
    \abs{x-x_0}^{\alpha(d,p)} & p<d, \\
     -\log \abs{x-x_0} & p=d,
  \end{cases} \qquad \mbox{ as } x\to x_0,
\end{equation}
where $f\sim g$ means that $$ \lim_{x\to x_0}\frac{f(x)}{g(x)}=
C$$ for some positive constant $C$.
\end{lemma}
\begin{remark}{\em
The asymptotics (\ref{nonremovasymp}) has been proved  for $p=2$ in
\cite{GiS}, for $1<p\leq d$ and $V=0$ in \cite[Theorem~2.1]{KV},
and also in some other cases in \cite{GV}. The proof below uses a
technique involving a scaling argument together with a comparison
principle that has been used for example in \cite{VBV}.}
\end{remark}
\begin{proof}
Assume that $1<p<d$, the proof for $p=d$ needs some minor
modifications, and is left to the reader. Without loss of
generality, we assume also that $x_0=0$.

Since $V\in L_{\mathrm{loc}}^\infty(\Omega)$, the solution $v$
satisfies $v(x)\asymp \abs{x}^{\alpha(d,p)}$. Let
 \be\label{maxseq}c:=\limsup_{x\to 0}\frac{v(x)}{\abs{x}^{\alpha(d,p)}}
=\lim_{n\to\infty}\frac{v(x_n)}{\abs{x_n}^{\alpha(d,p)}},\ee and
set
 $\mu(x):=c\abs{x}^{\alpha(d,p)}$. Define
$$v_n(x):=\abs{x_n}^{-\alpha(p,d)}v(\abs{x_n} x),$$
where $x_n\to 0$ is defined by \eqref{maxseq}.

Note that in an arbitrarily large punctured ball
$$C^{-1}\mu(x)\leq v_n(x)\leq C\mu(x)$$ for all $n$ large enough,
and in such a ball $v_n$ is a positive solution of the quasilinear
elliptic equation
$$-\Delta _p v_n(x) +\abs{x_n}^pV(x/\abs{x_n})v_n^{p-1}(x)=0.$$
Since $\{v_n\}$ is locally bounded and bounded away from zero in
any punctured ball, a standard elliptic argument implies that
there is a subsequence of $\{v_n\}$ that converges to a positive
singular solution $U$ of the limiting equation $-\Delta _p U=0$ in
the punctured space. Since $U\asymp \mu$ in the punctured space,
it follows that $U$ tends to zero at infinity. On the other hand,
\cite[Theorem~2.1]{KV}, implies that $U(x)\sim \mu(x)$ as $x \to
0$. Hence we can apply the comparison principle
(Theorem~\ref{CP}), and compare the functions $U$ and $\mu$ on
arbitrarily large balls, to obtain that $U= \mu$. This implies
that
\begin{equation}\label{eqvmu}
\lim_{n\to\infty} \| {v(x)/\mu (x)-1}\|_{L^\infty(\abs x=\abs{
x_n})}=0.\end{equation} In other words, $v$ is almost equal to
$\mu$ on a sequence of concentric spheres converging to $0$.

In order to prove that $v$ is almost equal to $\mu$ uniformly in
the sequence of the concentric annuli
 $A_n:=\{\abs{ x_n}\leq \abs{x}\leq \abs{
x_{n+1}}\}$, we construct two radial perturbations of $\mu$. Let
$\mu_-(x):=\mu(x)-\delta \abs x^{a}$ and $\mu_+(x):=\mu(x)+\delta
\abs x^{a}$ (for some $a>(p-d)/(p-1)$). It turns out that $\mu_-$
(resp., $\mu_+$) is a radial subsolution (resp., supersolution) of
the equation $Q'(u)=0$ near the origin, and therefore using the
comparison principle in the annulus $A_n$ and \eqref{eqvmu},  it
follows that
$$\lim_{r\to 0} \| {v(x)/\mu (x)-1}\|_{L^\infty(\abs x=\abs{
r})}=0.$$
\end{proof}

\begin{definition} {\em
Let $K$ be a compact set in $\Omega$.  A positive solution of the
equation $Q'(u)=0$ in $\Omega\setminus K$ is said to be a {\em
positive solution of minimal growth in a neighborhood of infinity
in} $\Omega$, if for any compact set $K_1$ in $\Omega$, with a
smooth boundary, satisfying $ \mathrm{int}(K_1)\supset K$, and any
positive supersolution $v\in C((\Omega\setminus K_1)\cup
\partial K_1)$ of the equation $Q'(u)=0$ in $\Omega\setminus K_1$,
the inequality $u\le v$ on $\partial K_1$ implies that $u\le v$ in
$\Omega\setminus K_1$. A positive solution $u$ of the equation
$Q'(u)=0$ in $\Omega$, which has minimal growth in a neighborhood
of infinity in $\Omega$ is called a {\em global minimal solution
of the equation $Q'(u)=0$ in $\Omega$}.}
\end{definition}

The following result is an extension to the $p$-Laplacian of the
corresponding result of S.~Agmon concerning  positive solutions of
real linear second-order elliptic operators \cite{Agmon}.

\begin{theorem}\label{thmmingr}
Suppose that $1<p\leq d$, and $Q$ is nonnegative on $\core$.
Then for any
$x_0\in \Omega$ the equation $Q'(u)=0$ has (up to a multiple
constant) a unique positive solution $v$ in
$\Omega\setminus\{x_0\}$ of minimal growth in a neighborhood of
infinity in $\Omega$.

 Moreover, $v$ is either a global minimal solution of the
equation $Q'(u)=0$ in $\Omega$, or $v$ has a nonremovable
singularity at $x_0$.
\end{theorem}

\begin{proof} Take $x_0\in \Omega$, and consider an exhaustion
$\{\Omega_{N}\}_{N=1}^{\infty}$ of $\Omega$ (as in the proof of
Theorem~\ref{pos}). Fix $N\geq 1$, and denote
$\Omega_{N,k}:=\Omega_N\setminus B(x_0,1/k)$.

Let  $\{f_k\}$ be a sequence of nonzero nonnegative smooth
functions such that for each $k\geq 2$, the function $f_k$ is
supported in $B(x_0,2/k) \setminus B(x_0,1/k)$.

Recall that $\lambda_{1,p}(\Omega_{N,k})>0$ for all $N,k\geq 1$.
By Theorem~\ref{thmGS}, there exists a unique positive solution of
the problem
\begin{equation}
\left\{
\begin{array}{rcl}
 Q'(u_{N,k})&=&c_kf_k \qquad \mbox{ in }
\Omega_{N,k},\\
u_{N,k}&=&0 \qquad\mbox{ on }\partial\Omega_{N,k},\\
u_{N,k}(x_1)&=&1,
\end{array}\right.
\end{equation}
where $x_1\neq x_0$ is a fixed point in $\Omega_1$ and $c_k>0$.
By Harnack's inequality and elliptic regularity, it follows that
$\{u_{N,k}\}$ admits a subsequence which converges locally
uniformly in $\Omega_N\setminus\{x_0\}$ to a positive solution
$G_N(\cdot,x_0)$ of the equation $Q'(u)=0$ in
$\Omega_N\setminus\{x_0\}$. Moreover, $G_N(\cdot,x_0)=0$ on
$\partial\Omega_{N}$, and $G_N(x_1,x_0)=1$.

Since $\lambda_{1,p}(\Omega_{N})>0$, and it is the unique
eigenvalue with a positive Dirichlet eigenfunction, it follows
that $G_N(\cdot,x_0)$ has a nonremovable singularity at $x_0$.

Recall that (\ref{nonremovasymp}) holds true with
$v(\cdot)=G_N(\cdot,x_0)$. It is convenient to normalize $G_N$ in
the traditional way, so that
 \be\label{nonremov2}\begin{array}{rccl}
      \displaystyle\lim_{x\to x_0}\frac{G_N(x,x_0)}{|x-x_0|^{\alpha(d,p)}}
      &=&\displaystyle\frac{p-1}{d-p}\left|S^{d-1}\right|^{-1/(p-1)} & \quad p<d,
    \\[6mm]
     \displaystyle\lim_{x\to x_0}\frac{G_N(x,x_0)}{-\log |x-x_0|}
     &=&\!\!\!\!\!\!\!\!\!\!\!\!\left|S^{d-1}\right|^{-1/(d-1)} & \quad
     p=d,
\end{array}
  \ee
  where $S^{d-1}$
  is the unit sphere in $\mathbb{R}^d$. Using the comparison principle
(Theorem~\ref{CP}), and (\ref{nonremovasymp}), it follows that
$G_N(\cdot,x_0)$ is the unique function with the above properties.
Hence, $G_N(\cdot,x_0)$ might be called \emph{the (Dirichlet)
positive $p$-Green function of the functional $Q$ in $\Omega_N$
with a pole at $x_0$}.

By the comparison principle (Theorem~\ref{CP}) and
(\ref{nonremovasymp}), it follows that the sequence
$\{G_N(\cdot,x_0)\}$ is nondecreasing as a function of $N$,  and
therefore it converges locally uniformly in
$\Omega\setminus\{x_0\}$ either to a positive function
$G(\cdot,x_0)$ or to infinity.

In the first case $G(\cdot,x_0)$ is a positive solution of the
equation $Q'(u)=0$ in $\Omega\setminus\{x_0\}$ and has the
asymptotic behavior (\ref{nonremovasymp}) near $x_0$. We call
$G(\cdot,x_0)$ \emph{the minimal positive $p$-Green function of
the functional $Q$ in $\Omega$ with a pole at $x_0$}.

In the second case, we consider the normalized sequence
$$v_N(x):=\frac{G_N(x,x_0)}{G_N(x_1,x_0)}\qquad n=1,2,\ldots.$$
By Harnack's inequality and elliptic regularity, it follows that
$\{v_{N}\}$ admits a subsequence which converges locally uniformly
to a positive solution $v$ of the equation $Q'(u)=0$ in
$\Omega\setminus\{x_0\}$.

Assume that $v$ has a nonremovable singularity at $x_0$.
Therefore, for each $N\geq 1$ we obtain by the comparison
principle and \eqref{nonremov2} that
$$G_N(x,x_0)\leq Cv(x)\qquad \forall x\in \Omega\setminus\{x_0\}$$
for some $C>0$ independent of $N$. But this contradicts our
assumption that $G_N\to \infty$ as $N\to \infty$.

Note that $G(\cdot,x_0)$ (in the first case) and $v$ (in the
second case) are limits of a sequence of positive solutions that
for $\delta>0$ are uniformly bounded on $\partial B(x_0,\delta)$,
and take zero boundary condition on $\partial \Omega_N$.
Therefore, by the comparison principle, $G$ and $v$ are positive
solutions in $\Omega\setminus\{x_0\}$ of minimal growth in a
neighborhood of infinity in $\Omega$. In particular, $v$ is a
global minimal positive solution of the equation $Q'(u)=0$ in
$\Omega$.

Using again the comparison principle and (\ref{nonremovasymp}), it
follows that such a solution is unique.
\end{proof}
The next theorem demonstrates that a global minimal positive
solution of the equation $Q'(u)=0$ in $\Omega$ is a ground state.

\begin{theorem}
Assume that $1<p\leq d$ and that $Q_V\geq0$ on $\core$. Then $Q_V$
is degenerately positive in $\Omega$ if and only if the equation
$Q'(u)=0$ admits a global minimal positive solution in $\Omega$.
\end{theorem}
\begin{proof}
Assume that $Q_V$ is strictly positive and assume that there
exists a global minimal positive solution $v$ of the equation
$Q'(u)=0$ in $\Omega$. By Proposition~\ref{strictpos}, there
exists a nonzero nonnegative function $V_1\in \core$ with $\supp
V_1\subset B(x_0,\delta)$ for some $\delta>0$, such that
$Q_{V-V_1}$ is strictly positive in $\Omega$. Therefore, in light
of Theorem~\ref{pos}, there exists a positive solution $v_1$ of
the equation $Q'_{V-V_1}(u)=0$ in $\Omega$.

Clearly, $v_1$ is a positive supersolution of the equation
$Q'_V(u)=0$ in $\Omega$ which is not a solution. On the other
hand, $v$ is a positive solution of  the equation $Q'_V(u)=0$ in
$\Omega$ which has minimal growth in a neighborhood of infinity in
$\Omega$. Therefore, there exists $\varepsilon>0$ such that
$\varepsilon v\leq v_1$ in $\Omega$. Define
$$\varepsilon_0:=\max\{\varepsilon>0\mid \varepsilon v\leq v_1 \mbox{ in } \Omega\}. $$
Clearly $\varepsilon_0 v\lneqq v_1$ in $\Omega$. Consequently,
there exist $\delta_1,\delta_2>0$ and $x_1\in \Omega$ such that
$$(1+\delta_1)\varepsilon_0 v(x)\leq v_1(x) \qquad x\in B(x_1,\delta_2).$$
Hence, by the definition of minimal growth, we have
$$(1+\delta_1)\varepsilon_0 v(x)\leq v_1(x) \qquad x\in \Omega\setminus B(x_1,\delta_2),$$
and thus $(1+\delta_1)\varepsilon_0 v\leq v_1$ in $\Omega$, which
is a contradiction to the definition of $\varepsilon_0$.

Assume that $Q$ admits a positive minimal $p$-Green function
$G(\cdot,x_0)$ in $\Omega$. We need to prove that $Q$ is strictly
positive.

Consider an exhaustion $\{\Omega_{N}\}_{N=1}^{\infty}$ of
$\Omega$ such that $x_0\in\Omega_1$ and $x_1\in
\Omega\setminus\Omega_1$. Fix a nonzero nonnegative function $f\in
C_0^\infty(\Omega_1)$. By Theorem~\ref{thmGS}, there exists a
unique positive solution of the Dirichlet problem
\begin{equation}
\left\{
\begin{array}{rcl}
 Q'(u_{N})&=&f \qquad \mbox{ in }
\Omega_{N},\\
u_{N}&=&0 \qquad\mbox{ on }\partial\Omega_{N}.
\end{array}\right.
\end{equation}
By the comparison principle (Theorem~\ref{CP}), $\{u_N\}$ is an
increasing sequence. Suppose that $\{u_N(x_1)\}$ is bounded. Then
$u_N\to u$, where $u$ satisfies the equation $Q'(u)=f\gvertneqq 0$
in $\Omega$. Since $u$ is a positive supersolution of the equation
$Q'(u)=0$ in $\Omega$ which is not a solution, Theorem~\ref{main}
implies that $Q$ is strictly positive.

Suppose that $u_N(x_1)\to\infty$. Then $v_N(x):=u_N(x)/u_N(x_1)$
solves the problem
\begin{equation}
\left\{
\begin{array}{rcl}
 Q'(v_{N})&=&\displaystyle\frac{f(x)}{u_N(x_1)^{p-1}} \qquad \mbox{ in }
\Omega_{N},\\[4mm]
v_{N}&=&0 \qquad\qquad\quad\mbox{ on }\partial\Omega_{N},\\[2mm]
 v_{N}(x_1)&=&1.
\end{array}\right.
\end{equation}
By Harnack's inequality, and the comparison principle,
\be\label{vngn}v_N \asymp G_N(\cdot,x_0) \qquad \mbox { in }
\Omega_N\setminus \Omega_1.
 \ee By a standard elliptic argument, we
may extract a subsequence of $\{v_N\}$ that converges to a
positive supersolution $v$ of the equation $Q'(u)=0$ in $\Omega$.
Recall that $G_N(\cdot,x_0)\to G(\cdot,x_0)$. Hence, (\ref{vngn})
implies that $v\asymp G(\cdot,x_0)$ in $\Omega\setminus \Omega_1$,
and in particular, $v$ is a positive solution in $\Omega\setminus
\Omega_1$ of minimal growth in a neighborhood of infinity in
$\Omega$. Note that $v\neq cG(\cdot,x_0)$ since $v$ is not
singular at $x_0$. Since the equation $Q'(u)=0$ does not admit a
global minimal solution in $\Omega$, it follows that $v$ satisfies
$Q'(u)\gvertneqq 0$ in $\Omega$, and by Theorem~\ref{main}, $Q$ is
strictly positive in $\Omega$.
\end{proof}
\mysection{Open problems}\label{sec:open} We conclude the paper
with a number of open problems suggested by the above results
which are left for future investigation. All these questions are
already resolved when $p=2$.

The first problem (Problem~\ref{pr:conv}) deals with the weakly
lower semicontinuity and convexity of the functional $Q$. As was
shown in Proposition~\ref{prwlsc}, $Q$ is weakly lower
semicontinuous in a limited sense, a property which is closely
related to convexity. However, $Q$ is not necessarily convex even
if $Q\ge 0$ on $\core$: for $p>2$ see the elementary
one-dimensional counterexample at the end of \cite{dPEM}, and also
the proof of Theorem~7 in \cite{GS}, for $p<2$ see \cite{FHTdT},
Example~2. Note also that Proposition~\ref{prwlsc} does not assert
that the domain $D$ defined by \eqref{defD} is weakly closed.
\begin{problem}\label{pr:conv}
 Does $Q$ have
a natural extension to a weakly closed set where it is weakly
lower semicontinuous? Under what conditions a nonnegative
functional $Q$ of the form \eqref{Q} is convex?
\end{problem}
Clearly $Q$ is convex when
$V\ge 0$ or $p=2$. The convexity of $Q$ gives rise to an energy space for the form
$Q$ that would generalize the space  $\mathcal{D}^{1,p}$,
similarly to the known case $p=2$ \cite{PT2}. Indeed, if $Q$ is a
nonnegative convex functional on $\core$, then it follows that
$Q^\frac{1}{p}$ is a norm on $\core$. Moreover, by
Theorem~\ref{main}, if $Q$ has a weighted spectral gap in
$\Omega$, then the completion of $\core$ with respect to this norm
is continuously imbedded into $L^p_\mathrm{loc}(\Omega)$. On the
other hand, if $Q$ has a ground state $v$ in $\Omega$, then $v$
belongs to the equivalence class of $0$ in this completion, and
there is no continuous imbedding of the completion even into
$\mathcal{D}^\prime(\Omega)$. However, due to (\ref{Poinc}), the
completion of $\core$ with respect to the norm induced by the
right hand side of (\ref{Poinc}) is continuously imbedded into
$L^p_\mathrm{loc}(\Omega)$.
\begin{problem}
Do the results of this paper extend to quasilinear functionals of
the form
$$Q^A(u):=\int_\Omega \left(|A(x)\nabla u|^{p}+V|u|^q\right)\dx,$$
where $A$ is a strictly positive definite matrix and $1<p\le q<
\infty$?
 \end{problem}
We note that Picone-type identity for the case $A(x)=a(x)I$ was
established in \cite{JTY}.
%

\begin{problem} Generalize  the results of Section~\ref{sec:mingr} to the case $d<p<\infty$.
 \end{problem}

\begin{problem}\label{pr4}
Let $\Omega_1\varsubsetneq\Omega$ be domains in $\R^d$. Suppose
that $Q$ is strictly positive in $\Omega_1$. Show that  there
exists an open domain $\Omega_1\subsetneqq\Omega_2\subset \Omega$
such that $Q_V$ is strictly positive in $\Omega_2$.
\end{problem}
Problem~\ref{pr4} was studied in \cite{PT1} under the assumption
$p=2$, and stronger statements were proved.

\begin{problem} Let $\Omega=\R^d$ and assume that $Q=Q_{V_j}$ are
strictly positive in $\Omega$, for $j=1,2$. For $y\in \R^d$ denote
$V_y(x):=V_1(x)+V_2(x-y)$. Show that under suitable decay
conditions on $V_j$, there exists $R>0$ such that for every
$y\in\R^d\setminus B_R(0)$ the functional $Q_{V_y}$ is strictly
positive in $\Omega$.
\end{problem}
This phenomenon has been proved for $p=2$ in \cite{KS,OS,S} for
Schr\"odinger operators and in \cite{Plocal} for the
non-selfadjoint case.

\appendix
\mysection{Appendix: Energy inequality}\label{sec:app} The
following inequality, established for $p\geq  2$, estimates the
functional $Q$ from below by an expression that leads to an
alternative proof of Lemma~\ref{someballs} for the case $p\geq 2$.

\begin{lemma}\label{Q12} Assume that $p\geq 2$. Let $v\in C^1(\Omega)$ be a positive solution
of the equation $Q'(u)=0$ in $\Omega$, and let $u\in\sobloc$,
$u\ge 0$, $\supp u\Subset\Omega$. Then
 \be\label{frombelow} Q(u)\;\; \left\{
                        \begin{array}{ll}
                          \geq Q_1(u)+Q_2(u) & \mbox{ if }\; p>2, \\[3mm]
                          =Q_1(u)  & \mbox{ if }\; p=2.
                        \end{array}
                      \right.
 \ee
where
 $$Q_1(u):=\frac 2p\int_\Omega |\nabla
v|^{p-2}v^2\left|\nabla\left[\left(\frac{u}{v}\right)^\frac{p}{2}\right]\right|^2\dx,\;
\mbox{ and} \quad Q_2(u):=\int_{\{\nabla v=0\}}|\nabla
u|^p\!\dx.$$
\end{lemma}

\begin{proof}
Since  $p\ge 2$, the obvious inequality
$t^p+(p-1)-pt\ge(p-1)(t-1)^2$ implies

\be \label{L1eq} L_1(u,v)\ge (p-1) \left(\frac{|\nabla
u|}{u}-\frac{|\nabla v|}{v}\right)^2\left(\frac{u}{v}\right)^p
|\nabla v|^{p-2}v^2,\ee where $L_1$ is defined by (\ref{L1}).
 We use the identity
 \be \label{mana}
\left|\nabla\left[\left(\frac{u}{v}\right)^\frac{p}{2}\right]\right|^2=
\left(\frac{p}{2}\right)^2\left(\frac{u}{v}\right)^p\left|\frac{\nabla
u}{u}-\frac{\nabla v}{v}\right|^2. \ee
Substitution of (\ref{mana}) into (\ref{L1eq}) and using the
identity \be \left(\frac{|\nabla u|}{u}-\frac{|\nabla
v|}{v}\right)^2-\left|\frac{\nabla u}{u}-\frac{\nabla
v}{v}\right|^2=2\frac{\nabla u\cdot\nabla v-|\nabla u||\nabla
v|}{uv} \ee
gives
\bean L_1(u,v)\ge (p-1)\left(\frac{2}{p}\right)^{2} |\nabla
v|^{p-2}v^2\left|\nabla\left[\left(\frac{u}{v}\right)^\frac{p}{2}\right]\right|^2
+\\[2mm] 2(p-1) (\nabla u\cdot\nabla v-|\nabla u||\nabla
v|)\left(\frac{u}{v}\right)^{p-1}|\nabla v|^{p-2},\eean
 which is the same as

\be\label{ina5} L_1(u,v)\ge (p-1)\left(\frac{2}{p}\right)^{2}
|\nabla
v|^{p-2}v^2\left|\nabla\left[\left(\frac{u}{v}\right)^\frac{p}{2}\right]\right|^2
- 2\frac{p-1}{p}L_2(u,v),\ee where $L_2$ is defined by (\ref{L2}).
Since $1\geq 2/(p^2-p)$ for $p\geq 2$, (\ref{ina5}) implies that

\be \label{2.13} L(u,v)\ge \frac{2}{p}|\nabla
v|^{p-2}v^2\left|\nabla\left[\left(\frac{u}{v}\right)^\frac{p}{2}\right]\right|^2.\ee
Note that for $p=2$ we have $L(u,v)= v^2\left|\nabla\left(\frac{u}{v}\right)\right|^2$ (see
\cite[Lemma 2.4]{PT2}).

Moreover, on the critical set $\{x\in\Omega\mid |\nabla
v(x)|=0\}$, we have  $L(u,v)=|\nabla u|^p$. Therefore for $p>2$ we have
$$L(u,v) \geq \frac2p|\nabla
v|^{p-2}v^2\left|\nabla\left[\left(\frac{u}{v}\right)^\frac{p}{2}\right]\right|^2+\mathbf{1}_{\{|\nabla
v|=0\}}|\nabla u|^p.$$ Integrating the latter inequality over
$\Omega$, we arrive at (\ref{frombelow}).
\end{proof}

\begin{remark}{\em
1. Note that (\ref{frombelow}) is based on the pointwise
inequality (\ref{2.13}).

2. Note that for $p=2$ we have $Q(u)=Q_1(u)$ (see for example
\cite[Lemma 2.4]{PT2}). The set of all
critical points of a positive solution of the linear equation $-\Delta
u+Vu=0$ in $\Omega$ has studied in \cite{HHHN}.}
\end{remark}

\noindent {\bf Alternative proof of Lemma~\ref{someballs}}.
 Assume that $p\geq 2$.  By (\ref{frombelow}), if $Q(u_k)\to 0$,
then $Q_1(u_k)\to 0$ and $Q_2(u_k)\to 0$.
Since $c_B=0$, there exists a sequence $u_k\in\core$, $u_k\ge 0$, such that
$\int_Bu_k^p=1$ and $Q(u_k)\to 0$. Repeating the first two steps
of the proof of Lemma~\ref{someballs}, we deduce that $\{u_k\}$ is
bounded in $W^{1,p}(\omega)$ for every $\omega\Subset\Omega$,
$\omega\supset B$.

Consider now a weakly convergent renamed subsequence
$u_k\rightharpoonup u$ in $W^{1,p}(\omega)$. Let
$$Q^\omega_1(w):=\!\frac2p\int_\omega\!\! |\nabla
v|^{p-2}v^2\left|\nabla\left[\left(\frac{w}{v}\right)^\frac{p}{2}\right]\right|^2\!\!\!\!\dx,
\; \mbox{ and }\;  Q^\omega_2(w):=\!\!\!\int_{\{\nabla
v=0\}\cap\omega}\!\!\!\!|\nabla w|^p\dx.$$

Since $Q^\omega_1$ and $Q^\omega_2$ are continuous convex
functionals on $W^{1,p}(\omega)$, they are weakly lower
semicontinuous, and therefore,
$$Q^\omega_1(u)\le\lim_{k\to\infty}
Q^\omega_1(u_k)=0\quad  \mbox{and} \quad
Q^\omega_2(u)\le\lim_{k\to\infty} Q^\omega_2(u_k)=0.$$

Consequently, $\nabla[\left(u/v\right)^\frac{p}{2}]=0$ almost
everywhere in $\omega\setminus\{\nabla v=0\}$ and $\nabla u=0$ in
$\omega\cap\{\nabla v=0\}$. Note that if $\nabla u(x)=0$ and
$\nabla v(x)=0$, then
$\nabla\left[\left(u(x)/v(x)\right)^\frac{p}{2}\right]=0$. Thus,
$\nabla\left[\left(u/v\right)^\frac{p}{2}\right]=0$ a.e. in
$\omega$, and since it holds for every $\omega$ containing $B$, it
follows that $u/v=\mathrm{constant}$ a.e. in $\omega$. By the
compact Sobolev imbedding on $B$,
$\int_Bu^p\dx=\lim\int_Bu_k^p\dx=1$, and therefore, $u=cv$, where
$c^{-p}=\int_B v^p\dx$. Note that any subsequence of $\{u_k\}$ has
a subsequence converging to $cv$ with the same $c$. From the
compactness of the local Sobolev imbedding, it follows that
$u_k\to cv$ in $L^p_{\mathrm{loc}}(\Omega)$. In other words,
$\{u_k\}$ is a null sequence. \qed

\begin{center}
{\bf Acknowledgments} \end{center} The authors are grateful to
L.~V\'{e}ron for kindly providing the authors the proof of
Lemma~\ref{lemveron}, a contribution which has substantially
improved the coherence of the paper. The authors wish to thank
V.~Liskevich and V.~Moroz for valuable discussions, in particular
for V.~Moroz' comment on the convexity counterexamples in
\cite{dPEM}, \cite{GS} and \cite{FHTdT}. Part of this research was done while
K.~T. was visiting the Technion and the Hebrew University, and
Y.~P. was visiting the Forschungsinsitut f\"{u}r Mathematik, ETH.
The authors would like to thank these institutes for the kind
hospitality. The work of Y.~P. was partially supported by the RTN
network ``Nonlinear Partial Differential Equations Describing
Front Propagation and Other Singular Phenomena",
HPRN-CT-2002-00274, and the Fund for the Promotion of Research at
the Technion. The work of K.~T. was partially supported by the
Swedish Research Council.

\end{document}